\documentclass[12pt]{article}
\usepackage[dvips]{graphicx}
\usepackage{amsfonts,amsmath,amssymb}
\usepackage{amsmath,amsthm}
\usepackage{latexsym}
\usepackage{color}
\usepackage{tikz}
\usepackage{color}
\newtheorem{thm}{Theorem}%[section]
\newtheorem{ob}[thm]{Observation}
\newtheorem{lem}[thm]{Lemma}%[section]
\newtheorem{conj}{Conjecture}
\newtheorem{cor}[thm]{Corollary}
\newtheorem{prop}[thm]{Proposition}
\newtheorem{ques}{Question}

\def\vertex(#1){\put(#1){\circle*{2}}}
\def\vertexo(#1){\put(#1){\circle{2}}}
\def\vert(#1){\put(#1){\circle*{1.5}}}
\def\verto(#1){\put(#1){\circle{1.5}}}
\def\lab(#1)#2{\put(#1){\makebox(0,0)[c]{#2}}}
\newcommand{\str}{\, \boxtimes \,}

\newcommand{\cp}{\,\Box\,}

\begin{document}

\title{On well-dominated direct, {Cartesian} and strong product graphs}

\author{Douglas F. Rall
\\ \\
Department of Mathematics \\
Furman University \\
Greenville, SC, USA\\
\small \tt Email: doug.rall@furman.edu}

\date{}
\maketitle
\begin{abstract}
If each minimal dominating set in a graph is a minimum dominating set, then the graph is called \emph{well-dominated}.  Since the seminal paper
on well-dominated graphs appeared in 1988, the structure of 
well-dominated graphs from several restricted classes have been studied.  In this paper we give a complete characterization of nontrivial
direct products that are well-dominated.  We prove that if a strong product is well-dominated, then both of its factors are well-dominated.
When one of the factors of a strong product is a complete graph, the other factor being well-dominated is also a sufficient condition for the 
product to be well-dominated.  Our main result gives a complete characterization of well-dominated Cartesian products in which at least one
of the factors is a complete graph.  In addition, we conjecture that this result is actually a complete characterization of the class of
nontrivial, well-dominated Cartesian products. 
\end{abstract}

{\small \textbf{Keywords:} well-dominated, Cartesian product, direct product, strong product } \\
\indent {\small \textbf{AMS subject classification:} 05C69, 05C76}

\section{Introduction} \label{sec:intro}
A dominating set of a finite graph $G$ is a set $D$ of vertices such that each vertex of $G$ is within distance $1$ of at least one vertex in $D$.
Finding a (set inclusion) minimal dominating set is straightforward and can be accomplished in linear time by simply ordering the vertices and
then discarding them one at a time if the remaining set still dominates the graph.  Depending on the graph, the resulting minimal dominating set may be
much larger than the domination number of $G$, which is the minimum cardinality among all its dominating sets.  On the other hand, for a given graph $G$
and positive integer $k$ it is well known that determining whether $G$ has a dominating set of cardinality at most $k$ is an $NP$-complete problem.
In most applications finding the size of a smallest dominating
set is the typical goal.  \emph{Well-dominated} graphs are those for which the above algorithm always finds a dominating set of minimum cardinality.
The seminal paper on well-dominated graphs was by Finbow, Hartnell and Nowakowski~\cite{fhn-1988}.  They characterized the well-dominated graphs of girth at
least $5$ and showed that the only well-dominated bipartite  graphs are those with domination number one-half their order.  Several other groups of authors
have studied the concept of well-dominated within restricted graph classes.  See, for example, \cite{tv-1990} (block and unicyclic graphs),
\cite{ptv-1996} (simplicial and chordal graphs), \cite{gks-2011} ($4$-connected, $4$-regular claw-free), and~\cite{fb-2015} (planar triangulations).
In~\cite{ghm-2017} G\"{o}z\"{u}pek, Hujdurovi\'{c} and Milani\v{c} characterized the well-dominated graphs that are nontrivial lexicographic products.

Our focus in this paper is the class of well-dominated graphs that have a nontrivial factorization as a Cartesian, direct or strong product.
These three graph products are referred to as the ``fundamental products'' in the book~\cite{hik-2011} by Hammack, Imrich and Klav\v{z}ar.  Along with the
lexicographic product they are the most studied graph products in the literature.

Anderson, Kuenzel and Rall~\cite{akr-2021}  characterized the direct products  that are well-dominated under the assumption that at least one of the factors has
no isolatable vertices.  (A vertex $x$ in a graph $X$ is \emph{isolatable} if there is an independent set $A$ in $X$ such that $\{x\}$ is a component
in $X-N[A]$.)

\begin{thm}~{\rm \cite[Theorem 3]{akr-2021}} \label{thm:directnoisolatables}
Let $G$ and $H$ be nontrivial connected graphs such that at least one of $G$ or $H$ has no isolatable vertices.  The
direct product $G \times H$ is well-dominated if and only if $G=H=K_3$ or at least one of the factors is $K_2$ and the other factor is
a $4$-cycle or the corona of a connected graph.
\end{thm}

In this paper we complete the characterization of well-dominated direct products by removing the requirement on isolatable vertices.  In fact, we
prove that if both factors of a direct product have order at least $3$ and one of them has an isolatable vertex, then the direct product is not
well-dominated.  The main result on well-dominated direct products is then the following theorem.

\begin{thm}  \label{thm:directcharacterization}
Let $G$ and $H$ be connected graphs.  The direct product $G \times H$ is well-dominated if and only if $G\times H=K_3\times K_3$,
$G\times H=K_2\times C_4$, or $G \times H=K_2 \times (F \odot K_1)$ for a connected graph $F$.
\end{thm}

In the same paper Anderson et al. proved that if a Cartesian product $G \cp H$ is well-dominated, then at least one of $G$ or $H$ is well-dominated.
In addition, they provided a characterization of the well-dominated Cartesian products of triangle-free graphs.  Namely, they proved that the Cartesian
product of two connected, triangle-free graphs of order at least $2$ is well-dominated if and only if both factors are complete graphs of order $2$.
Here we explore the more general case of well-dominated Cartesian products in which (at least) one of the factors has girth $3$.  In particular, we prove the following characterization of well-dominated graphs of the form $K_m \cp H$, and we conjecture that every well-dominated Cartesian product of nontrivial connected graphs is isomorphic to one of these.
\begin{thm} \label{thm:Cartesiancharacterization}
Let $m$ be a positive integer with $m \ge 2$ and let $H$ be a nontrivial, connected graph.  The Cartesian product $K_m \cp H$ is 
well-dominated if and only if either $m \neq 3$ and $H=K_m$ or $m=3$ and $H \in \{K_3,P_3\}$.
\end{thm}

For the strong product we prove that both factors of a well-dominated strong product are well-dominated.  If one of the factors of a strong product is a complete 
graph, then the other factor being well-dominated is also a sufficient condition for the product to be well-dominated.
\begin{thm} \label{thm:strongcompletefactor}
Let $n$ be a positive integer.  For any graph $H$ the strong product $K_n \str H$ is well-dominated if and only if
$H$ is well-dominated.
\end{thm}

The remainder of the paper is organized in the  following way.  In the next section we provide the necessary definitions for the remainder of the paper.  
In Section~\ref{sec:strong} we settle the relatively straightforward result for strong products.  Theorem~\ref{thm:directcharacterization},
the complete characterization of well-dominated direct products, is verified in Section~\ref{sec:direct}.  Proving Theorem~\ref{thm:Cartesiancharacterization} is the main task of Section~\ref{sec:Cartesian}.  We also derive a number of necessary conditions on two connected  graphs whose Cartesian product is well-dominated, which 
leads us to conjecture that the characterization in Theorem~\ref{thm:Cartesiancharacterization} captures all connected well-dominated Cartesian products.

\section{Definitions} \label{sec:defns}

All graphs in this paper are finite, undirected, simple and have order at least $2$.   For a positive
integer $n$, we let $[n]=\{1,\ldots,n\}$.  This set will be the vertex set of the complete graph of order $n$.
In general, we follow the terminology and notation of Hammack, Imrich, and Klav\v{z}ar~\cite{hik-2011}.  The order of a graph $G$ is the number of vertices in $G$
and is denoted $n(G)$; $G$ is \emph{nontrivial} if $n(G) \ge 2$.  For a vertex $v$ in a graph $G$, the \emph{open neighborhood}
$N(v)$ and the \emph{closed neighborhood} $N[v]$ are defined by $N(v)=\{u \in V(G)\,:\, uv \in E(G)\}$ and $N[v]=N(v)\cup \{v\}$.  For $A \subseteq V(G)$ we let
$N(A)=\cup_{v \in A}N(v)$ and $N[A]=N(A) \cup A$.  Any vertex subset $D$ such that $N[D]=V(G)$ is a \emph{dominating set} of $G$, and $D$ is then a \emph{minimal}
dominating set if no proper subset of $D$ is a dominating set.  The \emph{domination number} of $G$ is denoted by $\gamma(G)$ and is the minimum cardinality
among the dominating sets of $G$.  The \emph{upper domination number}, denoted $\Gamma(G)$, is the largest cardinality of a minimal dominating set of $G$.
A set $M \subseteq V(G)$ is an \emph{independent} set if its vertices are pairwise non-adjacent.  An independent set is  \emph{maximal} if it is not a proper subset
of an independent set.  The cardinalities of a smallest and a largest maximal independent set in $G$ are denoted by $i(G)$ and $\alpha(G)$, respectively.
Note that a maximal independent set is a dominating set, which gives
%\[\gamma(G) \le i(G) \le \alpha(G) \le \Gamma(G)\,.\]
\begin{equation} \label{eqn:domchain}
\gamma(G) \le i(G) \le \alpha(G) \le \Gamma(G)\,.
\end{equation}
A graph $G$ is called \emph{well-covered} if $i(G)=\alpha(G)$ and is \emph{well-dominated} if $\gamma(G)=\Gamma(G)$.  It is clear from~\eqref{eqn:domchain}
that the class of well-dominated graphs is a subclass of the class of well-covered graphs.

Let $G$ be a graph and let $u \in A \subseteq V(G)$.  The \emph{private neighborhood of $u$ with respect to $A$} is the set pn$[u,A]$ defined by
pn$[u,A]=\{x \in V(G):\,N[x]\cap A=\{u\}\}$.  Equivalently,  pn$[u,A]=N[u]-N[A-\{u\}]$.  The vertices in pn$[u,A]$ are called \emph{private neighbors}
of $u$ with respect to $A$.  The subset $A \subseteq V(G)$ is \emph{irredundant} if pn$[u,A]\neq \emptyset$ for every $u \in A$.
It follows from the definitions that a dominating set $D$ of $G$ is a minimal dominating set if and only if every vertex in $D$ has a private
neighbor with respect to $D$; that is, $D$ is irredundant.  In this paper we will also need a more restricted type of private neighbor.  The
\emph{external private neighborhood} of $u$ with respect to $A$ is the set epn$[u,A]$ defined by  epn$[u,A]={\rm pn}[u,A]-\{u\}=N(u)-N[A-\{u\}]$.
If $v \in {\rm epn}[u,A]$, then $v$ is called an \emph{external} private neighbor of $u$ with respect to $A$.  Such a vertex $v$, if it exists,
belongs to $V(G)-A$, which is the reason to use the word external.    A property related to irredundance, and one that is important for this paper,
is that of being open irredundant.  The set $A$ is \emph{open irredundant} if every vertex of $A$ has an external private neighbor with respect to $A$.

A vertex $x$ in a graph $G$ is an \emph{isolatable vertex} in $G$ if there exists an independent set $I$ of vertices in $G$ such that $x$ is isolated
in the induced subgraph $G-N[I]$ of $G$. Concerning the closed neighborhood of an independent set we have the following useful fact about well-dominated
graphs first observed by Finbow, Hartnell and Nowakowski.  The proof is straightforward and follows from the fact that if $M$ is an independent set in $G$
and $D$ is a minimal dominating set of $G-N[M]$, then $D \cup M$ is a minimal dominating set of $G$.

\begin{ob} {\rm \cite{fhn-1988}} \label{ob:hereditary}
If $G$ is a well-dominated graph and $M$ is any independent set of vertices in $G$, then $G-N[M]$ is well-dominated.
\end{ob}

Let $G$ and $H$ be finite, undirected graphs.  The \emph{Cartesian product} of $G$ and $H$, denoted $G \cp H$, has as its vertex set the Cartesian (set)
product $V(G) \times V(H)$. The \emph{direct product}, denoted $G \times H$, and the \emph{strong product}, $G \str H$, also have $V(G) \times V(H)$ as their
set of vertices.   Distinct vertices $(g_1,h_1)$ and $(g_2,h_2)$ are adjacent in
\begin{itemize}
\item $G \cp H$ if either ($g_1=g_2$ and $h_1h_2 \in E(H)$) or ($h_1=h_2$ and $g_1g_2 \in E(G)$);
\item $G \times H$ if $g_1g_2 \in E(G)$ and $h_1h_2 \in E(H)$;
\item $G \str H$ if they are adjacent in $G \cp H$ or they are adjacent in $G \times H$.
\end{itemize}

All three of these graph products are associative and commutative.  A product graph is called \emph{nontrivial} if both of its 
factors are nontrivial.  See~\cite{hik-2011} for specific information on these and other graph products.  The
\emph{corona} of a graph $G$, denoted by $G \odot K_1$, is the graph of order $2 n(G)$ obtained by adding, for each vertex $u$ of $G$ a new vertex $u'$
together with a new edge $uu'$.

\section{Well-dominated strong products} \label{sec:strong}

Recall that two vertices $(g_1,h_1)$ and $(g_2,h_2)$ are adjacent in the strong product $G \str H$ if one of the following holds.
\begin{itemize}
\item $g_1=g_2$ and $h_1h_2 \in E(H)$
\item $h_1=h_2$ and $g_1g_2 \in E(G)$
\item $g_1g_2 \in E(G)$ and $h_1h_2 \in E(H)$
\end{itemize}

Nowakowski and Rall~\cite{nr-1996} established the following relationships between ordinary domination invariants on strong products.
\begin{prop} {\rm \cite[Corollary 2.2]{nr-1996}} \label{prop:strongdomination}
If $G$ and $H$ are finite graphs, then
\[\gamma(G \str H) \le \gamma(G)\gamma(H) \text{ and } \Gamma(G \str H) \ge \Gamma(G)\Gamma(H)\,.\]
\end{prop}

The following corollary follows immediately from Proposition~\ref{prop:strongdomination}.
\begin{cor} \label{cor:strongw-cAndw-d}
Let $G$ and $H$ be finite graphs.  If $G \str H$ is well-dominated, then $G$ and $H$ are well-dominated.
\end{cor}

The converse of Corollary~\ref{cor:strongw-cAndw-d} is not true in general.  For example, the $5$-cycle is well-dominated, but
\[\gamma(C_5\str C_5)=4<6=\Gamma(C_5 \str C_5)\,.\]
However, we are able to show that if at least one of the factors is a complete graph, then the strong product of
this complete graph and a well-dominated graph is well-dominated.

Let $D$ be a dominating set of $G \str H$ and let $A =\{g \in V(G)\,:\,(g,h)\in D \text{ for some } h \in V(H)\}$.
If $u \in V(G)-A$ and $v \in V(H)$, then from the definition of
the edge structure of the strong product it follows that $(u,v)$ is dominated by $D$ only if $u$ is dominated by $A$.
Thus, $A$ dominates $G$ and therefore $\gamma(G) \le |A| \le |D|$.  Interchanging the roles of $G$ and $H$ proves the following lemma.

\begin{lem} \label{lem:strongprojection}
For all pairs of graphs $G$ and $H$, we have $\gamma(G \str H) \ge \max\{\gamma(G),\gamma(H)\}$.
\end{lem}

We now proceed to prove Theorem~\ref{thm:strongcompletefactor}, which is restated here.

\noindent \textbf{Theorem~\ref{thm:strongcompletefactor}} \emph{
Let $n$ be a positive integer.  For any graph $H$ the strong product $K_n \str H$ is well-dominated if and only if
$H$ is well-dominated.
}
\proof
Suppose that $H$ is a well-dominated graph.  Let $D$ be any minimal dominating set of $K_n \str H$
and let $S=\{h \in V(H)\,:\,(i,h)\in D \text{ for some } i \in [n]\}$.  Note that for any $u \in V(H)$ and for $1 \le i<j\le n$,
we have $N[(i,u)]=N[(j,u)]$.  Since $D$ is an irredundant set, we infer that $|S|=|D|$.
We claim that $S$ is a minimal dominating set of $H$.  As in the paragraph preceding Lemma~\ref{lem:strongprojection}
we see that $S$ dominates $H$.  Let $h$ be any vertex of $S$ and let $i \in [n]$ such that $(i,h)\in D$.  If $(i,h)$ is isolated in the subgraph
of $K_n \str H$ induced by $D$, then $h$ is isolated in the subgraph of $H$ induced by $S$ and $h \in {\rm pn}[h,S]$.
On the other hand, if $(i,h) \not\in {\rm pn}[(i,h),D]$, then there exists $(k,x) \in {\rm pn}[(i,h),D]$ such that $x \notin S$.
Again by definition of the edge structure of the strong product it
follows that $x \in {\rm pn}[h,S]$.  We conclude that $S$ is irredundant in $H$, and hence $S$ is a minimal dominating set of $H$.  Therefore,
$|D| = |S|=\gamma(H)$.  That is, all minimal dominating sets of $K_n \str H$ have the same cardinality, which implies that $K_n \str H$ is well-dominated.

The converse follows from Corollary~\ref{cor:strongw-cAndw-d}. \qed

Still unanswered is the following natural question.
\begin{ques}
What properties on the well-dominated graphs $G$ and $H$ are necessary and sufficient for $G \str H$ to be well-dominated?
\end{ques}

\section{Well-dominated direct products} \label{sec:direct}

In this section we complete the characterization of well-dominated direct products.  Throughout we assume that the factors are connected and have order at least $2$.
We prove Theorem~\ref{thm:directcharacterization}, which we now restate for convenience of the reader.

\noindent \textbf{Theorem~\ref{thm:directcharacterization}} \emph{
Let $G$ and $H$ be connected graphs.  The direct product $G \times H$ is well-dominated if and only if $G\times H=K_3\times K_3$,
$G\times H=K_2\times C_4$, or $G \times H=K_2 \times (F \odot K_1)$ for a connected graph $F$.
}

It is easy to see that if  $M$ is any maximal independent set of $G$ then $M \times V(H)$ is a maximal independent set, and thus also a minimal
dominating set, of $G \times H$.  If, in addition, $G \times H$ is well-dominated, then it follows that $M \times V(H)$ is a minimum dominating
set of $G \times H$ and so $\gamma(G \times H)=\alpha(G)n(H)$.

\begin{lem} \label{lem:directwell-dominated}
Let $G$ and $H$ be connected graphs of order at least $3$.  If $G \times H$ is well-dominated, then
$\gamma(G \times H)=\gamma(G)n(H)=\gamma(H)n(G)$. Furthermore, $\gamma(G)=\alpha(G)$ and $\gamma(H)=\alpha(H)$.
\end{lem}
\proof
Let $D$ be a minimum dominating set of $G$.  It is easy to see that $D \times V(H)$ dominates $G \times H$.    We get
\[\gamma(G)n(H)=|D \times V(H)|\ge \gamma(G \times H)=\alpha(G)n(H)\ge \gamma(G)n(H)\,,\]
and we have equality throughout.  Therefore, $\gamma(G \times H)=\gamma(G)n(H)$, and $\gamma(G)=\alpha(G)$.   By reversing the roles of $G$ and $H$
we also have $\gamma(G \times H)=\gamma(H)n(G)$ and  $\gamma(H)=\alpha(H)$.  \qed

As we have observed several times, $D \times V(H)$ dominates the direct product $G \times H$ if $D$ dominates $G$.  The next lemma shows that
if both factors have order at least $3$, then the set $D \times V(H)$ is not a \emph{minimal} dominating set unless $D$ is independent in $G$.

\begin{lem} \label{lem:allgammasetsindependent}
Let $G$ and $H$ be connected graphs of order at least $3$, and let $D$ be a minimal dominating set of $G$.  The set $D \times V(H)$ is a minimal
dominating set of $G \times H$ if and only if $D$ is independent in $G$.
\end{lem}
\proof
If $D$ is an independent dominating set of $G$, then $D \times V(H)$ is a maximal independent set, and hence a minimal dominating set, of $G \times H$.
For the converse, suppose that $a$ is a vertex that is not isolated in the subgraph of $G$ induced by $D$.  Let $x$ be a vertex in $H$ that is not a
support vertex.    We will show that ${\rm pn}[(a,x), D \times V(H)]=\emptyset$.
First, since $N(a) \cap D \neq \emptyset$, we see that $(a,x)$ has a neighbor in $D \times V(H)$, which implies that $(a,x) \not\in {\rm pn}[(a,x), D \times V(H)]$.
Next, let $(u,y) \in N((a,x))-(D \times V(H))$.  Since $x$ is not a support vertex of $H$, we see that there exists a path, say $x,y,z$, in $H$.
However, this means that $(a,z)$ is adjacent to $(u,y)$, and hence $(u,y) \not\in {\rm pn}[(a,x), D \times V(H)]$.  It follows that
${\rm pn}[(a,x), D \times V(H)]=\emptyset$, which implies that $D \times V(H)$ is not a minimal dominating set.  \qed

As a result of Lemma~\ref{lem:allgammasetsindependent} we can now state a very restrictive condition that must be satisfied by both
factors of a well-dominated direct product if both have order at least $3$.
\begin{cor} \label{cor:noedgesingammaset}
Let $G$ and $H$ be connected graphs of order at least $3$.  If $G \times H$ is well-dominated, then all minimum dominating sets of $G$
and all minimum dominating sets of $H$ are independent.
\end{cor}
\proof
Suppose $G$ and $H$ are connected of order at least $3$ such that $G \times H$ is well-dominated.  By Lemma~\ref{lem:directwell-dominated},
we have $\gamma(G \times H)=\gamma(G)n(H)=\gamma(H)n(G)$.  Let $D$ be any minimum dominating set of $G$.  Since $D \times V(H)$ dominates
$G \times H$ and $\gamma(G \times H)=|D \times V(H)|$, it follows from Lemma~\ref{lem:allgammasetsindependent} that $D$ is independent.
Similarly, every minimum dominating set of $H$ is independent. \qed

Using these results we now proceed to the proof of Theorem~\ref{thm:directcharacterization}.

It was verified in~\cite{akr-2021} that $G\times H$ is well-dominated if  $G\times H=K_3\times K_3$, $G\times H=K_2\times C_4$,
or $G \times H=K_2 \times (F \odot K_1)$ for a connected graph $F$.  For the converse we assume that $G \times H$ is well-dominated.
Suppose for the sake of contradiction that
$G\times H \neq K_3\times K_3$, $G\times H \neq K_2\times C_4$, and $G \times H \neq K_2 \times (F \odot K_1)$ for a connected graph $F$.
By Theorem~\ref{thm:directnoisolatables},  it follows that both $G$ and $H$ possess an isolatable vertex, and hence $n(G) \ge 3$ and $n(H) \ge 3$.
Let $M$ be an independent set in $G$ such that $G-N[M]=\{x\}$.  Since $G$ is connected, there exists a vertex $x'\in N(M)$ that is adjacent to $x$.
Let $S=M \cup \{x'\}$.  It is clear that $S$ dominates $G$.  Hence, $S$ is a minimum dominating set of $G$ since $|S|=|M|+1=\alpha(G)=\gamma(G)$
by Lemma~\ref{lem:directwell-dominated}, which is a contradiction by Corollary~\ref{cor:noedgesingammaset}.  \qed

Let $F$ be a connected graph.  It is straightforward to verify that $(F \odot K_1) \times K_2=(F \times K_2)\odot K_1$.  Thus the well-dominated graphs
of the form $K_2 \times (F \odot K_1)$ are a subclass of the easily recognizable well-dominated, bipartite coronas from the result of Finbow et al.~\cite{fhn-1988}
mentioned in Section~\ref{sec:intro}.

\section{Well-dominated Cartesian products} \label{sec:Cartesian}

In this section we prove the following characterization of well-dominated Cartesian products in which at least one of the factors 
is a complete graph.  In particular, we prove our main theorem in this study.

\noindent \textbf{Theorem~\ref{thm:Cartesiancharacterization}} \emph{
Let $m$ be a positive integer with $m \ge 2$ and let $H$ be a nontrivial, connected graph.  The Cartesian product $K_m \cp H$ is
well-dominated if and only if either $m \neq 3$ and $H=K_m$ or $m=3$ and $H \in \{K_3,P_3\}$.
}

We begin by deriving some preliminary results that will prove to be useful in its proof.  After that we prove Proposition~\ref{prop:completefactor},
which is a special case of Theorem~\ref{thm:Cartesiancharacterization} that assumes $n(H) \le m$. 

As mentioned in the introduction, the first step in this process was obtained by Anderson, Kuenzel and Rall~\cite{akr-2021}.

\begin{thm}~{\rm \cite[Theorem 1]{akr-2021}} \label{thm:beginningthm}
Let $G$ and $H$ be connected graphs.  If $G\cp H$ is well-dominated, then $G$ or $H$ is well-dominated.
\end{thm}

If a graph $G$ admits a dominating set that is open irredundant, then the following result provides a method for constructing
a minimal dominating set in any Cartesian product that has $G$ as one of its factors.

\begin{lem} \label{lem:minimalconstruction}
If  $D$ is an open irredundant dominating set of a graph $G$, then $D \times V(H)$
is an open irredundant, minimal dominating set of $G \cp H$, for any graph $H$.
\end{lem}
\proof  If $(x,y)\in V(G\cp H)-(D \times V(H))$, then there exists a vertex $x'\in N(x) \cap D$ since $D$ is a dominating set of $G$.  This
implies that $(x,y)$ is adjacent to  $(x',y)$, and thus $D \times V(H)$ is a dominating set of $G \cp H$.  To see that $D \times V(H)$
is a minimal dominating set, let $(d,h)$ be an arbitrary vertex in $D \times V(H)$.  Since $D$ is an open irredundant set of $G$,
there exists a vertex $d'$ in epn$[d,D]$.  By definition, $d'\in N(d)-D$, and $d'$ is not adjacent to any vertex of $D-\{d\}$.  Consequently,
$(d',h) \in {\rm epn}[(d,h), D \times V(H)]$.  Therefore, $D \times V(H)$ is an open irredundant, minimal dominating set of $G \cp H$.  \qed

Bollob\'{a}s and Cockayne proved that every graph with minimum degree at least $1$ has an open irredundant, minimum dominating set $D$.

\begin{prop} {\rm \cite[Proposition 6]{bc-1992}}\label{prop:openirredundant}
If a graph $G$ has no isolated vertices, then $G$ has a minimum dominating set that is open irredundant.
\end{prop}

By using Proposition~\ref{prop:openirredundant}, we now establish a relationship that must hold between two graphs if their Cartesian product is well-dominated.

\begin{prop} \label{prop:foundational}
Let $G$ and $H$ be nontrivial connected graphs.  If $G \cp H$ is well-dominated, then $\gamma(G \cp H)=\gamma(G)\cdot n(H)=\gamma(H)\cdot n(G)$.
\end{prop}
\proof  Let $D_1$ be a minimum dominating set of $G$ that is open irredundant and let $D_2$ be a minimum dominating set of $H$ that is open irredundant.
These minimum dominating sets of $G$ and $H$ exist by Proposition~\ref{prop:openirredundant}.  By Lemma~\ref{lem:minimalconstruction}, it follows that
both of $D_1 \times V(H)$  and $V(G) \times D_2$ are minimal dominating sets of $G \cp H$.  Since $G \cp H$ is well-dominated, it follows that
$\gamma(G \cp H)= |D_1 \times V(H)|=|V(G) \times D_2|$, and therefore
\[ \gamma(G \cp H)= \gamma(G)\cdot n(H)= \gamma(H)\cdot n(G)\,.\] \qed

If $D$ is a minimal dominating set of a graph $G$, then, in general, some vertices of $D$ will have external private neighbors with respect to $D$
and some vertices will not. The following notation will be useful in what follows.  If $D$ is any dominating set of $G$, we define $c(D)$ as follows:
\[c(D)=\{x \in D\,:\, {\rm pn}[x,D]=\{x\}\}\,.\]
That is, a vertex $x$ in $D$ belongs to $c(D)$ if and only if $x$ is isolated in the subgraph induced by $D$, and $N(x)-D \subseteq N(D-\{x\})$.
In particular, $c(D)$ is a (possibly empty) independent set of $G$.  Note also that if $c(D) \neq \emptyset$, then $|D| \geq 2$ since we
are assuming that $G$ is of order at least $2$.

\begin{lem}\label{lem:c(D)construction}
Let $S$ be a minimal dominating set of a graph $H$.  If $D$ is a minimal dominating set of a graph
$G$ such that $c(D) \neq \emptyset$, then $(\{u\} \times S) \cup \left( (D-\{u\}) \times V(H)\right)$ is a dominating set of $G \cp H$ for every $u \in c(D)$.
If, in addition, $c(D)=\{u\}$, then $(\{u\} \times S) \cup \left( (D-\{u\}) \times V(H)\right)$ is a minimal dominating set of $G \cp H$.
\end{lem}
\proof  Let $u \in c(D)$.  For simplification, let $A=(\{u\} \times S) \cup \left( (D-\{u\}) \times V(H)\right)$.  Suppose $(x,y)\in V(G\cp H)-A$.  If $x \not\in D$,
then there exists $d \in D-\{u\}$ such that $dx \in E(G)$.  This follows since $D$ dominates $G$ and  $x\notin {\rm pn}[u,D]=\{u\}$.  Hence, $(x,y)$
has a neighbor in $(D-\{u\}) \times V(H)$.  On the other hand, if
$x \in D$, then $x=u$ and $y \notin S$.  Since $S$ dominates $H$, it follows that $(x,y)$ has a neighbor in $\{u\} \times S$.  Therefore, $A$ is a dominating set
of $G \cp H$.

Now, suppose that $c(D)=\{u\}$.
 To show that $A$ is a minimal dominating set let $(a,b) \in A$.  If $a \neq u$, then $a$ has an external  private neighbor, say $a'$, with respect to $D$, and
$(a',b)$ is a  private neighbor of $(a,b)$ with respect to $A$.  On the other hand, suppose $a=u$.  This implies that $b \in S$.  Since $S$ is a minimal dominating
set of $H$, there exists $b' \in {\rm pn}[b,S]$.  It now follows that $(u,b')$ is a private neighbor of $(u,b)$ with respect to $A$ since $u$ is an isolated vertex
in $G[D]$.  We have shown that pn$[(a,b),A] \neq \emptyset$, and it follows that $A$ is a minimal dominating set of $G \cp H$.   \qed

If $G$ is any finite graph, then for large enough $m$, namely for $m>\Delta(G)$, the Cartesian product $G \cp K_m$ is well-covered.  
(See page 1262 of~\cite{hrw-2018}.)  This is not true in the well-dominated class as we now prove.

\begin{prop} \label{prop:NoWell-dominated}
If $G$ is a nontrivial connected graph and has a minimum dominating set $D$ such that $c(D)\neq \emptyset$, then
$G \cp H$ is not well-dominated for every nontrivial connected graph $H$.
\end{prop}
\proof Let $G$  be a nontrivial, connected graph and suppose that $D$ is a minimum dominating set of $G$ with a vertex
$u \in c(D)$.  Let $H$ be a nontrivial connected graph and suppose that $S$ is any  minimum dominating set of $H$.  By Lemma~\ref{lem:c(D)construction},
$(\{u\} \times S) \cup \left( (D-\{u\}) \times V(H)\right)$ is a  dominating set of $G \cp H$. While $(\{u\} \times S) \cup \left( (D-\{u\}) \times V(H)\right)$
may not be a minimal dominating set, it contains one.  Since
\begin{align*}
|(\{u\} \times S) \cup \left( (D-\{u\}) \times V(H)\right)| &= \gamma(H)+(\gamma(G)-1)\cdot n(H)\\
                                                            & = \gamma(G)\cdot n(H)+(\gamma(H)-n(H))< \gamma(G)\cdot n(H),
\end{align*}
we conclude by Proposition~\ref{prop:foundational} that $G \cp H$ is not well-dominated. \qed

The following corollary of Proposition~\ref{prop:NoWell-dominated} further limits which graphs can be a factor of a well-dominated Cartesian product.

\begin{cor} \label{cor:noisolatables}
If $G$ is a nontrivial connected graph that is well-dominated and has an isolatable vertex, then $G \cp H$ is not well-dominated for any
nontrivial connected graph $H$.
\end{cor}
\proof Suppose $G$ is a nontrivial connected, well-dominated  graph and suppose $x$ is an isolatable vertex of $G$.  Let $H$ be any connected graph of
order at least $2$.  Since $G$ is well-dominated, we have $\gamma(G)=i(G)=\alpha(G)=\Gamma(G)$.  Let $I$ be an independent set in $G$ such that
$G-N[I]=\{x\}$.  The set $J=I \cup \{x\}$ is an independent dominating set, and is therefore also a minimum dominating set of $G$.  Since $G-N[I]=\{x\}$,
we see that $x \in c(J)$.  It follows by Proposition~\ref{prop:NoWell-dominated} that $G \cp H$ is not well-dominated.  \qed

The following result follows immediately from Corollary~\ref{cor:noisolatables}.

\begin{cor} \label{cor:noleaves}
If $G$ is a connected, well-dominated graph of order at least $3$ and $G \cp H$ is well-dominated for some nontrivial connected graph $H$, then
$\delta(G) \ge 2$.
\end{cor}

Proceeding with the proof of Theorem~\ref{thm:Cartesiancharacterization} we first deal with the case where the order of
the complete factor in the Cartesian product is at least as large as the order of the other factor.

\begin{prop} \label{prop:completefactor}
Let $m$ be a positive integer larger than $1$ and let $H$ be a nontrivial connected graph such that $n(H) \le m$.  The Cartesian product $K_m \cp H$
is well-dominated if and only if one of the following holds.
\begin{enumerate}
\item $m=2$ and $H=K_2$.
\item $m=3$ and $H \in \{K_3,P_3\}$.
\item $m \ge 4$ and $H=K_m$.
\end{enumerate}
\end{prop}
\proof  If $n(H) < m$, then $K_m \cp H$ is not well-dominated follows immediately from Proposition~\ref{prop:foundational}.  Hence, we now assume that $H$ has order $m$.
It is straightforward to show the result is correct for $m=2$ or $m=3$.  Now, let $m \ge 4$.  It is easy to show that  $K_m \cp K_m$ is well-dominated.  For
the converse,  suppose that $H$ is a connected graph of order $m$ such that $K_m \cp H$ is well-dominated but such that $H \neq K_m$.
Throughout the proof we let $\{h_1,h_2,\ldots,h_m\}$ denote the vertex set of $H$.

By Proposition~\ref{prop:foundational}, we infer that $\gamma(H)=1$.  Without loss of generality we assume that $\{h_1\}$ dominates $H$.
Let $A$ be a maximum independent set of $H$.  Since $H$ is not a complete graph, $|A| \ge 2$ and $A$ does not contain $h_1$.
  If $|A|=m-1$, then $H=K_{1,m-1}$.
In this case let $S=([m-1]\times \{h_m\}) \cup (\{m\} \times \{h_2,\ldots,h_{m-1}\})$.  Note that $|S|=2m-3>m$.  We claim that
$S$ is a minimal dominating set of $K_m \cp H$.  It is clear that $S$ dominates $K_m \cp H$.    Furthermore,
pn$[(i,h_m),S]=\{(i,h_1)\}$ for $i\in [m-1]$, and  $(m,h_j)\in {\rm pn}[(m,h_j),S]$ for $2 \le j \le m-1$.
This proves that $S$ is a minimal dominating set of $K_m\cp H$, which is a contradiction.
Thus, we may assume that $|A|<m-1$, and we may also assume that $A=V(H)-\{h_1,\ldots,h_k\}$ for some $k$ such that $2 \le k \le m-2$.  Let $R=([m-1] \times \{h_1\}) \cup (\{m\} \times A)$.  It is easy to show that $R$ dominates $K_m \cp H$.  For $i \in [m-1]$, we see that $(i,h_2) \in {\rm pn}[(i,h_1),R]$, and
$(m,a)\in {\rm pn}[(m,a),R]$, for  $a \in A$.    This implies that $R$ is a minimal dominating set.  However, $|R|=m-1+|A|\ge m-1+2>m$, which
is a contradiction and implies that $K_m \cp H$ is not well-dominated.

Therefore, if $m \ge 4$ and $H$ is a connected graph of order $m$, then $K_m \cp H$ is well-dominated if and only if $H=K_m$.   \qed

In the remainder of the proof of Theorem~\ref{thm:Cartesiancharacterization} we consider $K_m \cp H$ where $n(H) > m$.  The proofs for $m=2$ and $m=3$
are straightforward, while the proof of the more general case for $m \ge 4$ occupies the rest of this section.

\begin{prop} \label{prop:connectedm=2}
If $H$ is a connected graph of order at least $3$, then $K_2 \cp H$ is not well-dominated.
\end{prop}
\proof  Suppose there exists a connected graph $H$ of order at least $3$, such that $K_2 \cp H$ is well-dominated.  By Proposition~\ref{prop:foundational}, we
have $\gamma(K_2 \cp H)=\gamma(K_2)n(H)=n(H)$.   Since $H$ has order at least $3$ and is connected, there exists a vertex $h \in V(H)$ such that
$\deg(h) \ge 2$.  Now, if $S$ is the set defined by  $S=([2] \times \{h\}) \cup (\{1\} \times (V(H)-N_H[h]))$ we arrive at a contradiction since
$S$ dominates $K_2 \cp H$ and $|S|\le n(H)-1$.  \qed

We have a similar result for Cartesian products with $K_3$.

\begin{prop} \label{prop:connectedm=3}
If $H$ is a connected graph of order at least $4$, then $K_3 \cp H$ is not well-dominated.
\end{prop}
\proof  Suppose there exists a connected graph $H$ of order more than $3$, such that $K_3 \cp H$ is well-dominated.  By Proposition~\ref{prop:foundational}, we
have $\gamma(K_3 \cp H)=\Gamma(K_3\cp H)=n(H)$.  Suppose first that $\Delta(H) \ge 3$; let $h$ be a vertex of $H$ such that $\deg(h)=r \ge 3$.  This implies
that $S=([3] \times \{h\}) \cup (\{1\} \times (V(H)-N_H[h]))$ dominates $K_3 \cp H$.  This is a contradiction since $|S|=3+(n(H)-(r+1)) < n(H)$.  Consequently,
$\Delta(H)=2$, and thus $H$ is either a cycle or a path of order at least $4$.  It is easy to verify that $\gamma(K_3 \cp P_4)=4 < 6=\Gamma(K_3 \cp P_4)$ and
$\gamma(K_3\cp C_4)=3<6=\Gamma(K_3 \cp C_4)$.  Thus, we may assume that $H$ contains a path of order $5$, say $h_1h_2h_3h_4h_5$.  Let
$A=(\{2,3\} \times \{h_3\}) \cup (\{1\} \times (V(H)-\{h_2,h_3,h_4\}))$.  This set $A$ dominates $K_3 \cp H$ and yet $|A|=n(H)-1$.  This final contradiction
establishes the proposition.  \qed

\begin{prop} \label{prop:gammalarger1}
Let $m$ be a positive integer such that $m \ge 4$. If $H$ is a connected graph of order more than $m$, then $K_m \cp H$ is not well-dominated.
\end{prop}
\proof  We proceed by induction on $m$.  For the base case we suppose for the sake of contradiction that there exists a connected graph $H$
of order more than $4$ such that $K_4 \cp H$ is well-dominated.  Let $k=\gamma(H)$.  By Proposition~\ref{prop:foundational},
$\gamma(K_4 \cp H)=n(H)=4k$, and hence every minimal dominating set of the well-dominated graph $K_4 \cp H$ has cardinality $4k$.
 Let $D=\{h_1,\ldots,h_k\}$ be a minimum dominating set of $H$.
Suppose first that $H$ has a vertex $x$ of degree at least $4$.  If $A=V(H)-N[x]$, then the set $S$ defined by $S=([4] \times \{x\}) \cup (\{1\} \times A)$
is a dominating set of $K_4 \cp H$.  However,  $|S|=4+|A|=4+(n(H)-|N[x]|) < 4k$, which is a contradiction.
Therefore, $\Delta(H)\le 3$.  Since $V(H)=\cup_{i=1}^kN[h_i]$, it follows that $\deg(h_i)=3$, for every $i \in [k]$, and we see that
$N_H[h_1], \ldots,N_H[h_k]$ is a partition of $V(H)$.  Let $M=\{4\} \times D$ and for each $i \in [k]$, let $X_i=N(h_i)=\{x_{i1},x_{i2},x_{i3}\}$ and let $Y_i=X_i \cup \{h_i\}$.
The set $M$  is independent in the well-dominated graph $K_4 \cp H$ and by Observation~\ref{ob:hereditary} it follows that the graph $G$ defined by
$G=K_4 \cp H - N[M]$ is well-dominated.  Note that $G = K_3 \cp F$, where $F$ is the subgraph of $H$ induced by $V(H)-D$.  This implies that each component of
$G$ is well-dominated.  Using Proposition~\ref{prop:connectedm=3} we infer that each component of $F$ has order at most $3$. Furthermore, $\Delta(F) \le 2$ since
$\Delta(H)=3$.

Note that $K_3 \cp K_2$ is not well-dominated, which then implies that each component of $F$ has order $1$ or $3$.  At least one of the components of $F$ has
order $3$, for otherwise $H$ is not connected.  Each component of $F$ that has order $3$ intersects either one, two or three of the sets $X_1,\ldots,X_k$.
Suppose first that there exists $i \in [k]$, say $i=1$, such that $\langle X_i\rangle$ is a component of $F$.  This implies that the subgraph of $H$ induced
by $X_1 \cup \{h_1\}$ is a component of $H$ of order $4$, which contradicts the assumption that $H$ is connected and has order at least $5$.  We thus assume
that each component of $F$ that has order $3$ has a nonempty intersection with either $2$ or $3$ of the sets $X_1,\ldots,X_k$.  Any such component is clearly
either a path of order $3$ or a complete graph of order $3$.

Suppose there exists $1 \le i<j\le k$ such that $\langle X_i \cup X_j \rangle$ contains a $P_3$ or $K_3$ involving at least one vertex from each
of $X_i$ and $X_j$.  Without loss of generality we assume that $i=1$ and $j=2$ and that $x_{13}$ has degree at least $2$ in $\langle X_1 \cup X_2 \rangle$.
We assume without loss of generality that $x_{13}x_{21} \in E(H)$.  If $x_{13}$ is adjacent to another vertex in $X_2$, say $x_{13}x_{22} \in E(H)$, then
let $B=( \{1\} \times \{h_1,x_{11}, x_{12},x_{23}\}) \cup \{(2,x_{21}),(3,x_{21}),(4,x_{22})\}$.  On the other hand, if $x_{13}$ is adjacent to another vertex in $X_1$, say $x_{13}x_{12} \in E(H)$, then let $B= (\{1\} \times \{h_1,x_{11}, x_{22}\}) \cup \{(2,x_{12}),(2,x_{23}),(3,x_{21}),(4,x_{21})\}$.  In both cases we
see that $B \cup ([4] \times \{h_3,\ldots,h_k\})$ dominates $K_4 \cp H$ and has cardinality $4k-1$, which is a contradiction.

Hence, every component of $F$ that has order $3$ contains one vertex from three distinct members of the partition $X_1,\ldots,X_k$ of $V(F)$.  We assume without loss
of generality that $x_{11},x_{21},x_{31}$ is a path in $F$.  Let
$B=(\{1\} \times \{h_1,h_2,x_{12},x_{13},x_{22}, x_{23},x_{32}\}) \cup \{(2,x_{11}),(2,x_{33}),(3,x_{31}),(4,x_{31})\}$.  It now follows that
$B \cup ([4] \times \{h_4,\ldots,h_k\})$ is a dominating set of $K_4 \cp H$ and has cardinality $4k-1$, which is a contradiction.
Therefore, if $H$ is a connected graph of order more than $4$, then $K_4 \cp H$ is not well-dominated.

Now let $m \ge 5$ and suppose that if $G$ is any connected graph of order at least $m$, then $K_{m-1} \cp G$ is not  well-dominated.  Again, for the sake of arriving at a contradiction, suppose there exists a connected graph $H$ of order more than $m$ such that $K_m\cp H$ is well-dominated.  For consistency and ease of understanding we use the same notation as in the case $m=4$.   Let $k=\gamma(H)$ and let $D=\{h_1,\ldots,h_k\}$ be a minimum dominating set of $H$.  By Proposition~\ref{prop:foundational},  $\gamma(K_m \cp H)=n(H)=mk$, and hence every minimal dominating set of the well-dominated graph $K_m \cp H$ has cardinality $mk$.  Suppose first that $H$ has a vertex $x$ of degree at least $m$.  If $A=V(H)-N[x]$, then the set $S$ defined by
$S=([m] \times \{x\}) \cup (\{1\} \times A)$  is a dominating set of $K_m \cp H$.  However,  $|S|=m+|A|=m+(n(H)-|N[x]|) < mk$, which is a contradiction.

Therefore, $\Delta(H)\le m-1$.  Since $V(H)=\cup_{i=1}^kN[h_i]$, it follows that $\deg(h_i)=m-1$, for every $i \in [k]$  and $N_H[h_1], \ldots,N_H[h_k]$ is a partition of $V(H)$.
Similar to the case above (for $m=4$) we let $M$ be the independent set defined by $M=\{m\} \times D$, and we note by Observation~\ref{ob:hereditary} that  $G=K_m\cp H - N[M]$  is well-dominated, and hence every component of $G$ is well-dominated.  Since $G$ is isomorphic to  $K_{m-1} \cp F$, where $F$ is the subgraph of $H$ induced by $V(H)-D$, it follows from the inductive hypothesis and Proposition~\ref{prop:completefactor} that every nontrivial component of $F$ is isomorphic to $K_{m-1}$.  Also, by Proposition~\ref{prop:foundational} we see that $\gamma(F)=k$, which implies that $F$ has no components of order $1$.  That is, $F$ is the disjoint union of $k$ complete graphs
of order $m-1$.  Suppose that for some $i \in [k]$ there exists a vertex $x \in N_H(h_i)$ such that $N_F(x) \subseteq N_H(h_i)$.  It follows that $N_H(h_i)$ induces  a
complete graph in $F$, which implies that $H$ is not connected.  This contradiction means there is no such $i \in [k]$.  In particular, for each $u\in N_H(h_1)$ we have
$N(u) \cap (N_H(h_2) \cup \cdots \cup N_H(h_k)) \neq \emptyset$.  Let
\[S=\left ( \{1\} \times \cup_{i=2}^kN_H[h_i]\right ) \cup (\{2, \ldots, m\} \times \{h_1\})\,.\]
This set $S$ dominates $K_m \cp H$ and $|S|=m(k-1) +(m-1)=mk-1$, which is a contradiction. This establishes the proposition.   \qed

By combining the results of Propositions~\ref{prop:completefactor},~\ref{prop:connectedm=2},~\ref{prop:connectedm=3}, and~\ref{prop:gammalarger1}, 
the proof of Theorem~\ref{thm:Cartesiancharacterization} is complete.

We close this section on well-dominated Cartesian products with the following conjecture.

\begin{conj} \label{conj:onefactorcomplete}
Let $G$ and $H$ be nontrivial connected graphs.  If $G \cp H$ is well-dominated, then at least one of $G$ or $H$ is a complete
graph.
\end{conj}

If Conjecture~\ref{conj:onefactorcomplete} is true, then we would have a complete characterization of the well-dominated Cartesian products.
That is, if Conjecture~\ref{conj:onefactorcomplete} is true, then by Theorem~\ref{thm:Cartesiancharacterization} it follows
that the Cartesian product $G\cp H$ of two nontrivial, connected graphs is well-dominated
if and only if $G \cp H=K_m \cp K_m$ for some positive integer $m\ge 2$ or $G \cp H=K_3 \cp P_3$.

\vfill
\end{document}